\documentclass{article}

\usepackage{amssymb}
\usepackage{amsfonts}
\usepackage{latexsym}

\newtheorem{prop}{Proposition}
\newtheorem{cor}{Corollary}

\newenvironment{pf}{{\bf Proof.}}{\hfill$\Box$\\[1mm]}

\newcommand{\A}{\mathcal{A}}

\renewcommand{\H}{\mathcal{H}}

\newcommand{\Sp}{\mathcal{S}}
\newcommand{\T}{\mathcal{T}}

\begin{document}
\begin{center}
{\LARGE Dirac Operators on\\
Quantum Flag Manifolds}\\[1cm]
{\large Ulrich Kr\"ahmer}\\[3mm]
Fakult\"at f\"ur Mathematik und Informatik\\
Universit\"at Leipzig, Augustusplatz 10, 04109 Leipzig, Germany\\
E-mail: kraehmer@math.uni-leipzig.de
\end{center}

\abstract{A Dirac operator $D$ on
quantized irreducible generalized flag manifolds is defined.
This yields a Hilbert space realization of the covariant first-order 
differential calculi constructed by 
I. Heckenberger and S. Kolb. All differentials 
$df=i[D,f]$ are bounded operators. 
In the simplest case of Podle\'s' quantum sphere 
one obtains the spectral triple found by L. Dabrowski and
A. Sitarz.}

\section{Introduction}
After about 20 years of quantum groups and of noncommutative
geometry in the sense of A. Connes the relation between these two theories is
still not understood very well.
In particular, there is no theory linking Connes' concept of 
spectral triple to that of finite-dimensional covariant differential calculus on 
quantum spaces as developed by S.L. Woronowicz [W]. It is only known 
that the basic examples of such calculi which are the
$3D$-calculus and the $4D_\pm$-calculi on the quantum group 
$SU_q(2)$ itself can \emph{not} be realized by spectral triples. This
was shown by K. Schm\"udgen in [S]. So the spectral triples 
of [CP, G] are not related to these calculi.\\
The aim of the present paper is to show that 
$q$-deformations of irreducible generalized flag manifolds 
$M=G/P$ behave better in this respect.\\
In [HK1] I. Heckenberger and S. Kolb proved that these quantum spaces admit 
exactly two irreducible finite-dimensional covariant differential calculi 
$(\Gamma_\pm,d_\pm)$. Their direct sum
$(\Gamma,d)$ is a $\ast$-calculus whose elements are
$q$-analogues of complex-valued differential forms on the real manifold $M$.\\
The main result of this paper is that $\Gamma$ can be realized by bounded
operators on a Hilbert space such that $df=i[D,f]$ for a self-adjoint
operator $D$. The latter generalizes the classical Dirac operator on $M$. 
A calculation of its spectrum seems 
to be a non-trivial problem, but its construction suggests that the spectrum is 
a smooth deformation of the classical one.\\ 
The simplest example of a generalized flag manifold is the
complex projective line $\mathbb{C}P^1 \simeq S^2$. The corresponding
quantum flag manifold is Podle\'s' quantum sphere in the so-called
quantum subgroup case. For this one obtains the Dirac operator found by
L. Dabrowski and A. Sitarz [DS1].\\
The paper is organized as follows: The first two sections are devoted
to background material on quantum groups and quantum flag
manifolds. In Sections~\ref{t},\ref{c} and \ref{s} we define
analogues of the tangent space, its 
Clifford algebra and the spinor bundle of $M$. With these as ingredients 
the Dirac operator $D$ is constructed in Section~\ref{d}. In the final
section we study the associated differential calculus and
prove the main results.\\[2mm] 
While working on this paper the author profited a lot from 
discussions with many people. He would like to
thank especially I. Agricola, T. Friedrich,
I. Heckenberger and S. Kolb for their suggestions and help.

\section{Quantum groups}\label{p}
In this section we fix notations and recall some definitions
and results of quantum group theory used in the sequel. We refer to
the monographs [KS] and [J] for proofs and further details.\\
Throughout this paper $\mathfrak{g}$ is a complex
simple Lie algebra, $G$ the corresponding 
connected simply connected Lie group, $\mathfrak{h}$ is a
Cartan subalgebra of $\mathfrak{g}$ and 
$\{\alpha_1,\ldots,\alpha_N\}$ and $\{\omega_1,\ldots,\omega_N\}$ are a
set of simple roots and the corresponding set of 
fundamental weights. The integral root and weight lattices are denoted 
by $\mathbf{Q}$ and $\mathbf{P}$, respectively. The dominant integral
weights are denoted by $\mathbf{P}^+$. The Killing form induces a
bilinear pairing $\langle \cdot , \cdot \rangle$ on 
$\mathbf{P} \times \mathbf{P}$. Then 
$\langle \omega_i,\alpha_j \rangle=:\delta_{ij} d_i$.\\
Let $U_q(\mathfrak{g})$ be the quantized
universal enveloping algebra corresponding to 
$\mathfrak{g}$ in the form denoted by $\check U$ in [J, 3.2.10]. We denote its
generators by $K_\lambda,E_i,F_i$, $i=1,\ldots,N$,
$\lambda \in \mathbf{P}$ and set
$K_i:=K_{\alpha_i}$. See [KS, Section~6.1.2] for their explicit
relations. For the coproduct, counit and antipode 
we use the conventions of [J]. In particular, the coproduct of
$E_i,F_i$ is 
$$
	\Delta(E_i)=E_i \otimes 1 + K_i \otimes E_i,\quad
	\Delta(F_i)=F_i \otimes K_i^{-1} + 1 \otimes F_i.
$$
We assume $q \in (1,\infty)$ and consider
$U_q(\mathfrak{g})$ as the Hopf $\ast$-algebra called 
the compact real form of $U_q(\mathfrak{g})$ in [KS, Section~6.1.7]. 
Its involution $\ast$ coincides 
on generators with $\kappa'$ from [J, 3.3.3], 
$$
	E_i^\ast=K_iF_i,\quad
	F_i^\ast=E_iK_i^{-1},\quad
	K_i^\ast=K_i,
$$
but $\kappa'$ is continued to a linear map. 
As in [GZ] we set $\theta := \ast \circ S$.\\
There is a bilinear form $\langle \cdot,\cdot \rangle$ on
$U_q(\mathfrak{g})$ called the quantum Killing form 
(or Rosso form) which is invariant under the adjoint action 
$$
	\mathrm{ad}(X)Y := X \triangleright Y := 
	X_{(1)} Y S(X_{(2)}),\quad X,Y \in U_q(\mathfrak{g})
$$
of $U_q(\mathfrak{g})$ on itself in the sense that 
$$
	\langle Z \triangleright X,Y \rangle =
	\langle X,S(Z) \triangleright Y \rangle \quad
	\forall X,Y,Z \in U_q(\mathfrak{g}).
$$
Above as in the following we use Sweedler notation 
$\Delta(X)=X_{(1)} \otimes X_{(2)}$ for the coproduct $\Delta$ of a Hopf
algebra.\\
The quantum Killing form satisfies [J, 3.3.3, 7.2.4]
$$
	\langle X^\ast,Y \rangle =
	\overline{\langle Y^\ast,X \rangle},\quad
	\langle X^\ast,X \rangle > 0\quad
	\forall X,Y \in U_q(\mathfrak{g}),X \neq 0.
$$
and vanishes on 
$U^\lambda_q(\mathfrak{g}) \times U^\mu_q(\mathfrak{g})$ if
$\lambda \neq -\mu$. Here
$U^\mu_q(\mathfrak{g})$, $\mu \in \mathbf{Q}$, consists of the elements 
$X \in U_q(\mathfrak{g})$ with
$K_\lambda X K_\lambda^{-1}=q^{\langle \lambda,\mu \rangle}X$ for all
$\lambda \in \mathbf{P}$.\\
The representation theories of $U_q(\mathfrak{g})$ and $\mathfrak{g}$ are
closely related. In particular, for every $\lambda \in \mathbf{P}^+$ there
exists an irreducible representation $(\rho_\lambda,V_\lambda)$ 
of highest weight $\lambda$ and a Hermitian
inner product $(\cdot,\cdot)_\lambda$ on $V_\lambda$ 
such that $\rho_\lambda$ becomes
a $\ast$-representation. The weight structure of $V_\lambda$ and the
decomposition of tensor products into irreducible
components remains unchanged as well [KS, Chapter~7].\\  
Let $\mathbb{C}_q[G]$ be the coordinate algebra of the standard 
quantum group associated to $G$. This is the Hopf $\ast$-subalgebra of
the Hopf dual $U_q(\mathfrak{g})^\circ$ generated by all
matrix coefficients $t_{ij}^\lambda$ of the representations 
$V_\lambda$, $\lambda \in \mathbf{P}^+$. We will conversely
treat elements of $U_q(\mathfrak{g})$ as functionals on 
$\mathbb{C}_q[G]$ and write the dual pairing between 
$X \in U_q(\mathfrak{g})$ and $f \in \mathbb{C}_q[G]$ as $X(f)$.\\
This pairing turns $\mathbb{C}_q[G]$ into a 
$U_q(\mathfrak{g})$-bimodule with left and
right action given by $X \triangleright f:=X(f_{(2)})f_{(1)}$,
$f \triangleleft X:=X(f_{(1)})f_{(2)}$. The structure of this bimodule is
given by the classical Peter-Weyl theorem. That is, the
$t_{ij}^\lambda$ form a vector space basis 
of $\mathbb{C}_q[G]$ and for fixed $i$ 
(fixed $j$) a basis of the representation $V_\lambda$ with respect to the
left (right) action.\\
The linear functional $h : \mathbb{C}_q[G] \rightarrow \mathbb{C}$
defined by $h(t^\lambda_{ij}):=\delta_{\lambda 0}$ and called the Haar
functional is biinvariant with respect to the
$U_q(\mathfrak{g})$-actions [KS, Section~11.3]. The associated
Hermitian inner product
$$
	\langle f,g \rangle_h := h(f g^\ast),\quad f,g \in \mathbb{C}_q[G]
$$
is positive definite and is the direct sum of the
$(\cdot,\cdot)_\lambda$ arising from considering $\mathbb{C}_q[G]$ as
right $U_q(\mathfrak{g})$-module. That is, we have
$$
	\langle f \triangleleft X,g \rangle_h = 
	\langle f,g	 \triangleleft X^\ast \rangle_h \quad
	\forall X \in U_q(\mathfrak{g}),f,g \in \mathbb{C}_q[G].
$$  
 
\section{Quantum flag manifolds}\label{f}
Let $P$ be a standard parabolic subgroup of $G$ and $M=G/P$ the
corresponding generalized flag manifold [FH, \S 23.3, BE].
As a real manifold $M$ is diffeomorphic to $G_0/L_0$, where
$G_0$ denotes the compact real form of $G$, $L$ is the Levi factor
of $P$ and $L_0:=L \cap G_0$, cf. [BE, 6.4].
Let $\mathfrak{p},\mathfrak{l}$ denote the Lie algebras of $P,L$. Throughout this
paper we assume that $M$ is irreducible, that is, that
$\mathfrak{g}/\mathfrak{p}$ is irreducible with respect to the adjoint
action of $\mathfrak{p}$. This implies that
$\mathfrak{l}$ is the Lie subalgebra of
$\mathfrak{g}$ generated by $\mathfrak{h}$ and the root vectors
$E_i,F_i$ associated to the simple roots $\alpha_i$, 
$i \neq r$, for a certain $r$, see [BE, Example~3.1.10]. For example, if
$G=SL(N+1,\mathbb{C})$, then the irreducible flag manifolds exhaust the complex 
Grassmann manifolds $Gr(r,N+1)$, $r=1,\ldots,N$.\\  
Let $U_q(\mathfrak{l}) \subset U_q(\mathfrak{g})$
be the Hopf $\ast$-subalgebra generated by 
$K_\lambda,E_i,F_i$, $\lambda \in \mathbf{P}$, $i \neq r$ and define
$$
	\mathbb{C}_q[M]:=
	\{f \in \mathbb{C}_q[G] \>|\> X \triangleright f = \varepsilon(X)f \>\forall
	X \in U_q(\mathfrak{l})\}.
$$
This algebra is a $q$-deformation of
a $\ast$-algebra of complex-valued functions on $M$.
Completion with respect to the $C^\ast$-completion of 
$\mathbb{C}_q[G]$ leads to a $q$-deformation
of the $C^\ast$-algebra associated to the
compact topological space $M$. With slight abuse of terminology from algebraic
geometry we call $\mathbb{C}_q[M]$ the coordinate algebra of the quantum flag
manifold $M_q$. See for example [HK1, DS2] for more information
about quantum flag manifolds.\\
The right and the left action of $U_q(\mathfrak{g})$
on $\mathbb{C}_q[G]$ commute. Hence 
$\mathbb{C}_q[M]$ is a right $U_q(\mathfrak{g})$-module. It
decomposes into irreducible components in the same way as its classical
analogue.
 
\section{The tangent space}\label{t}
Let $\mathfrak{u}$ be the orthogonal complement of $\mathfrak{l}$ 
with respect to the Killing form of $\mathfrak{g}$. It decomposes as 
$\mathfrak{u}=\mathfrak{u}_+ \oplus \mathfrak{u}_-$, where 
$\mathfrak{p}=\mathfrak{l} \oplus \mathfrak{u}_+$ is the Levi
decomposition of $\mathfrak{p}$ and $\mathfrak{u}_-$ can be identified
with the complex tangent space $\mathfrak{g}/\mathfrak{p}$ of 
$G/P$ at $eP$. The adjoint action of $\mathfrak{l}$ on $\mathfrak{u}$ 
defines an embedding of $\mathfrak{l}$ into 
$\mathfrak{so}(2m,\mathbb{C})$, where 
$\mathrm{dim}_\mathbb{C}\,M=m$. We now introduce 
analogues of $\mathfrak{u},\mathfrak{u}_\pm$  
for quantum flag manifolds.\\
Let $\lambda=-2n \cdot \omega_r$. The number 
$n \in \mathbb{N} \setminus \{0\}$ is arbitrary but fixed and will
play no role in the sequel. But probably it may be used to adjust
the analytical properties of the Dirac operator we will derive below.\\
Define $X_0:=K_\lambda-1$ and
$$
	X_1:=F_r \triangleright X_0 = F_r \triangleright K_\lambda=
	F_r K_\lambda K_r -
	K_\lambda F_r K_r=
	(1-q^{2nd_r}) F_r K_r K_\lambda. 
$$
\begin{prop}
The adjoint action turns
$\mathfrak{u}_- := \mathrm{ad}(U_q(\mathfrak{l})) X_1$
into the irreducible finite-dimensional representation of 
$U_q(\mathfrak{l})$ with highest weight $-\alpha_r$.
\end{prop}
\begin{pf}
Since $\Delta(K_\mu)=K_\mu \otimes K_\mu$ and
$S(K_\mu)=K_\mu^{-1}$ we have
$$
	K_\mu \triangleright X_1 =
	K_\mu X_1 K_\mu^{-1} =
	q^{-\langle \mu,\alpha_r \rangle} X_1 \quad
	\forall \mu \in \mathbf{P}.
$$
Furthermore, for $i \neq r$ we have
$$
	E_i \triangleright X_1
	= E_iF_r \triangleright K_\lambda
	= F_rE_i \triangleright K_\lambda
	= 0,
$$ 
because $K_\lambda$ commutes with all $E_i,F_i$ for $i \neq r$ and
therefore
$$
	E_i \triangleright K_\lambda =
	E_i K_\lambda -
	K_i K_\lambda K_i^{-1} E_i=
	E_i K_\lambda-E_i K_\lambda =
	0. 
$$
Since $X_0$ belongs to the locally finite
part of $U_q(\mathfrak{g})$ [J, 7.1.3] the vector space 
$\mathfrak{u}_-$ is finite-dimensional.
Hence the claim follows. 
\end{pf}
Fix a basis $X_i$ of $\mathfrak{u}_-$ consisting of weight vectors 
and define 
$X^i:=X_i^\ast$. Then the $X^i$ form a basis of a vector space
which we denote by $\mathfrak{u}_+$. Since
\begin{equation}\label{adast}
	(X \triangleright Y)^\ast=\theta(X) \triangleright Y^\ast 
\end{equation}
this vector space is $\mathrm{ad}(U_q(\mathfrak{l}))$-invariant as
well. Set $\mathfrak{u}:=\mathfrak{u}_+ \oplus \mathfrak{u}_-$.
We also introduce
$\mathfrak{u}_0:=\{X \in \mathfrak{u}\>|\>X^\ast=X\}$ as
an analogue of the real tangent space of $M$. It is 
invariant under $U_q(\mathfrak{l}_0):=\{X \in
U_q(\mathfrak{l})\>|\>\theta(X)=X\}$. Note that $U_q(\mathfrak{l}_0)$ is a
subalgebra of $U_q(\mathfrak{l})$, but not a Hopf subalgebra [GZ].\\
Since the highest weight representations of quantized universal enveloping
algebras are for $q$ not a root of unity of the same structure as their
classical counterparts the complex dimension of 
$\mathfrak{u}_\pm$ equals $m$.\\
The weights of $X_i$ are all distinct, so after an appropriate
normalization we have $\langle X_i,X^j \rangle = \delta_{ij}$.

\section{The Clifford algebra}\label{c}
We now define a quantum Clifford algebra associated to $M_q$.
We refer to [Fr] for the appearing notions from classical spin geometry.\\
The Clifford algebra $\mathrm{Cl}(2m,\mathbb{C})$ 
is the universal algebra for which there is a
vector space embedding 
$\gamma : \mathbb{C}^{2m}\rightarrow \mathrm{Cl}(2m,\mathbb{C})$
such that
$\gamma(v)^2=-\sum_{i} v_i^2$ for all $v \in \mathbb{C}^{2m}$. 
The spin representation $\sigma$ on the space 
$\Sigma_{2m}:=\mathbb{C}^{2^m}$ of $2m$-spinors yields an
isomorphism 
$\mathrm{Cl}(2m,\mathbb{C}) \simeq \mathrm{End}(\Sigma_{2m})$.\\ 
The crucial point leading to a quantum analogue of 
$\mathrm{Cl}(2m,\mathbb{C})$ in the context of quantum flag
manifolds is that $\gamma$ is 
$\mathfrak{so}(2m,\mathbb{C})$-equivariant. Here the
representation of $\mathfrak{so}(2m,\mathbb{C})$ on $\mathbb{C}^{2m}$ is the
vector representation $\rho$ and the one on 
$\mathrm{End}(\Sigma_{2m}) \simeq \Sigma_{2m} \otimes \Sigma_{2m}^\ast$ is 
the tensor product $\sigma \otimes \sigma^\ast$ 
of the spin representation and the
dual representation. In fact, the standard vector space isomorphism 
$\mathrm{Cl}(2m,\mathbb{C}) \simeq \Lambda^\ast \mathbb{C}^{2m}$ is an
isomorphism of $\mathfrak{so}(2m,\mathbb{C})$-representations and
$\gamma$ is the restriction to 
$\mathbb{C}^{2m}=\Lambda^1 \mathbb{C}^{2m}$.\\
Not all flag manifolds are spin manifolds [CG], but all admit 
spin$^\mathbb{C}$ structures [Fr, Section~3.4].
In any case the embedding 
$\mathfrak{l} \subset \mathfrak{so}(2m,\mathbb{C})$ defines the 
representations $\rho$ and $\sigma$ of $\mathfrak{l}$
and $\rho$ appears in $\sigma \otimes \sigma^\ast$.\\
These representations of $\mathfrak{l}$ can be deformed to representations of
$U_q(\mathfrak{l})$ which we denote by the same symbols. 
The decomposition of $\sigma \otimes \sigma^\ast$
into irreducible components remains the same. Hence we have:
\begin{prop}
There is a $U_q(\mathfrak{l})$-equivariant embedding 
$$
	\gamma : \mathfrak{u}_+ \oplus \mathfrak{u}_- \rightarrow
	\mathrm{End}(\Sigma_{2m}).
$$
\end{prop}
Without loss of generality we can assume that 
$$
	\gamma(X^i) = \overline{\gamma(X_i)}^T =:
	\gamma(X_i)^\ast, 
$$
because we could embed first only $\mathfrak{u}_-$ and take 
the above formula as the definition of $\gamma(X^i)$.\\
Note that the map $\gamma$ is not uniquely determined by these
conditions, but it always can be assumed to be a
smooth deformation of the classical one.\\
We call the algebra generated by $\gamma(\mathfrak{u}_0)$ and 
$\gamma(\mathfrak{u})$ the real and complex quantum Clifford algebra
associated to the quantum flag manifold $M_q$.

\section{The spinor bundle}\label{s}
Next we define a spinor bundle $\Sp$ over $M_q$ in form of a quantum 
homogeneous vector bundle \cite{GZ} generalizing the 
homogeneous vector bundle $G_0 \times_{L_0} \Sigma_{2m}$ over $M$.
It is defined in terms of the following vector space whose elements
are interpreted as its sections:
\begin{eqnarray}\label{spb}
	\Gamma(M_q,\Sp)
&:=& \{\psi \in \mathbb{C}_q[G] \otimes \Sigma_{2m}\>|\>
	X \triangleright \psi = \sigma(S(X))\psi \>\forall X \in 
	U_q(\mathfrak{l})\}\nonumber\\
&\simeq& \bigoplus_{\lambda \in \mathbf{P}^+}
	V_\lambda \otimes 
	\mathrm{Hom}_{U_q(\mathfrak{l})}(V_\lambda,\Sigma_{2m}).
\end{eqnarray} 
The isomorphism $\simeq$ is given by
Peter-Weyl decomposition of $\mathbb{C}_q[G]$.\\
If $\{A^\lambda_i\}$ are bases of 
$\mathrm{Hom}_{U_q(\mathfrak{l})}
(V_\lambda,\Sigma_{2m})$ for all $\lambda$ for
which this space is non-trivial and if the matrix
coefficients $t^\lambda_{ij}$ of the Peter-Weyl basis are
defined with respect to the bases $\{v^\lambda_j\}$ of $V_\lambda$,
then the elements
$$
	\psi^\lambda_{ij} := 
	\sum_{k} S(t^\lambda_{kj}) \otimes A^\lambda_i(v^\lambda_k)
$$
form a basis of $\Gamma(M_q,\Sp)$.\\ 
We define a Hermitian inner product 
$\langle \cdot , \cdot \rangle_\Sp$
on $\Gamma(M_q,\Sp)$ by applying 
$\langle \cdot,\cdot \rangle_h$ to $\mathbb{C}_q[G]$ and
the invariant Hermitian inner product 
$( \cdot,\cdot )_\sigma$ to $\Sigma_{2m}$.
We can choose the basis $\psi^\lambda_{ij}$ to be orthonormal. We
complete $\Gamma(M_q,\Sp)$ to a Hilbert space $\H$ which we call
the space of square-integrable spinor fields on the quantum flag
manifold $M_q$.\\
The quantized universal enveloping algebra $U_q(\mathfrak{g})$ acts on
$\Gamma(M_q,\Sp)$ from the right 
(by acting from the right on $\mathbb{C}_q[G]$). The multiplication in
$\mathbb{C}_q[G]$ defines a $\mathbb{C}_q[M]$-bimodule structure on 
$\Gamma(M_q,\Sp)$ and when restricting to a one-sided action one
obtains a projective module over $M_q$ [GZ]. 

\section{The Dirac operator}\label{d}
Let $D_-$ be the linear operator acting 
on $\mathrm{Hom}_{U_q(\mathfrak{l})}(V_\lambda,\Sigma_{2m})$
by
$$
	D_- : A \mapsto
 	- \sum_{i} \gamma(X^i) \circ A \circ \rho_\lambda(X_i).
$$
The following proposition shows that $D_-$ is well-defined.
\begin{prop}
We have $D_-(A) \in 
\mathrm{Hom}_{U_q(\mathfrak{l})}(V_\lambda,\Sigma_{2m})$.
\end{prop}  
\begin{pf}
For $Y \in U_q(\mathfrak{l})$ we have
\begin{eqnarray}
&& \sum_{i} \gamma(X^i) \circ A \circ \rho_\lambda(X_i)
	\rho_\lambda(S(Y))\nonumber\\
&=& \sum_{i} \gamma(X^i) \circ A \circ 
	\rho_\lambda(S(Y_{(1)})Y_{(2)} X_i S(Y_{(3)}))\nonumber\\
&=& \sum_{i} \gamma(X^i) \sigma(S(Y_{(1)})) \circ  
	A \circ \rho_\lambda(Y_{(2)} \triangleright X_i)\nonumber\\
&=& \sum_{ij} \gamma(X^i) \sigma(S(Y_{(1)})) \circ
	A \circ \rho_\lambda(\langle Y_{(2)} \triangleright X_i,X^j \rangle 
	X_j)\nonumber\\
&=& \sum_{ij} \gamma(
	\langle X_i,S(Y_{(2)}) \triangleright X^j \rangle
	X^i) \sigma(S(Y_{(1)})) \circ A \circ \rho_\lambda(X_j)\nonumber\\
&=& \sum_{j} \gamma(S(Y_{(2)}) \triangleright X^j)
	\sigma(S(Y_{(1)})) \circ A \circ \rho_\lambda(X_j)\nonumber\\
&=& \sum_{j} \sigma(S(Y_{(3)}))\gamma(X^j)\sigma(S^2(Y_{(2)}))
	\sigma(S(Y_{(1)})) \circ A \circ \rho_\lambda(X_j)\nonumber\\
&=& \sigma(S(Y)) \sum_{j} \gamma(X^j) \circ A \circ \rho_\lambda(X_j),\nonumber
\end{eqnarray} 
where we used 
the Hopf algebra axioms and the equivariance of $\gamma$.
\end{pf}
The resulting operator on
$\Gamma(M_q,\Sp)$ which acts trivially on $V_\lambda$ in (\ref{spb})
will be denoted by the same symbol. It can be extended to a linear operator on  
$\mathbb{C}_q[G] \otimes \Sigma_{2m}$ acting by
$$
	D_- : f \otimes v \mapsto 
	- \sum_{j} (S^{-1}(X_j) \triangleright f) \otimes 
	\gamma(X^j) v.
$$
We consider $D_-$ as densely defined operator on $\H$.
Analogously there is an operator $D_+$ acting on
$\mathbb{C}_q[G] \otimes \Sigma_{2m}$ by
$$
	D_+ : f \otimes v \mapsto 
	- \sum_{j} (S^{-1}(X^j) \triangleright f) \otimes 
	\gamma(X_j) v.
$$
Finally we define the Dirac operator $D:=D_++D_-$. Notice that for
$X \in U_q(\mathfrak{g})$ and $f,g \in \mathbb{C}_q[G]$ the
$U_q(\mathfrak{g})$-invariance of $h$ and (\ref{adast}) imply
\begin{eqnarray}
	h((X \triangleright f)g^\ast)
&=& h((X_{(1)} \triangleright f)(X_{(2)}S(X_{(3)}) \triangleright g^\ast))\nonumber\\
&=& h(X_{(1)} \triangleright (f (S(X_{(2)}) \triangleright
	g^\ast)))\nonumber\\
&=& h(f(S^2(X)^\ast \triangleright g)^\ast).\nonumber
\end{eqnarray} 
Hence $D$ is symmetric on the domain
$\Gamma(M_q,\Sp)$. It is the direct sum of its
restrictions to the spaces 
$V_\lambda \otimes 
\mathrm{Hom}_{U_q(\mathfrak{l})}(V_\lambda,\Sigma_{2m})$ which are all
finite-dimensional and pairwise orthogonal. Hence it
becomes diagonal in a suitable orthonormal basis and therefore extends
to a self-adjoint operator on $\H$ which we denote by $D$ as well.\\
It seems to be a non-trivial task to generalize Parthasarathy's
formula for $D^2$ [Pa] to the quantum case and to calculate explicitly
the spectrum of $D$ as in [CFG]. But since only finite matrices are
involved which are smooth deformations of those describing the
classical Dirac operator, the spectrum should as well be 
a smooth deformation of the classical spectrum.\\
In the simplest example of a generalized flag manifold
$M=\mathbb{C}P^1 \simeq S^2$ the corresponding quantum
flag manifold is one of Podle\'s' quantum spheres [Po]. In this case
L. Dabrowski and A. Sitarz derived a Dirac operator starting with an
ansatz and implementing the axioms for real spectral triples [DS1]. 
It follows from the uniqueness result [DS1, Lemma~5] that the Hilbert
space representation of $U_q(\mathfrak{sl}(2,\mathbb{C}))$ 
and $\mathbb{C}_q[M]$ calculated in [DS1] is the one on 
$\H$ considered here. Inserting the explicit formulae 
for the left action of $U_q(\mathfrak{sl}(2,\mathbb{C}))$ on
the Peter-Weyl basis one sees that the Dirac operators also coincide.

\section{The differential calculus}\label{k}
In this section we study the covariant first-order 
differential calculi over $\mathbb{C}_q[M]$ induced
by $D_\pm$ and $D$. We refer to [KS] and
[HK2] for the general theory of covariant differential calculi on quantum
groups and quantum homogeneous spaces.\\ 
The author's main impetus to study quantum flag manifolds was
the result of [HK1] mentioned already in the introduction 
that on quantum flag manifolds as discussed here there exist 
exactly two finite-dimensional irreducible (first-order) covariant differential
calculi $(\Gamma_\pm,d_\pm)$. 
These calculi have dimension $m$ and their direct
sum $(\Gamma,d)$ is a $\ast$-calculus.\\
A spectral triple $(\A,\H,D)$ over a $\ast$-algebra $\A$ (cf.~[C])
always induces a differential $\ast$-calculus with $df:=i[D,f]$, $f \in \A$.
The result of this section will be that the operators $D_\pm,D$
realize the calculi $\Gamma_\pm,\Gamma$ in this way by
bounded operators on $\H$.\\
We will treat only $\Gamma_-$, the analogous results for 
$\Gamma_+$ and $\Gamma$ are immediate.\\ 
In the rest of this paper, elements of $\mathbb{C}_q[M]$ will always
be treated as linear operators on $\H$ by considering the \emph{right} action of 
$\mathbb{C}_q[M]$ on $\Gamma(M_q,\Sp)$. If one rewrites this paper 
starting with the left coset space $P \setminus G$, the
constructions will work for the left action instead.\\
We denote by $\Gamma_-'$ 
the differential calculus over $\mathbb{C}_q[M]$ defined by $D_-$:
$$
	\Gamma_-':=\bigl\{i\sum_{j} f_j [D_-,g_j] 
	\>\bigr|\> f_j,g_j \in \mathbb{C}_q[M]\bigr\} \subset
	\mathrm{End}(\Gamma(M_q,\Sp)).
$$
Then the following formula for $d'_- f:=i[D_-,f]$ holds:
\begin{prop}\label{df}
For all $f \in \mathbb{C}_q[M]$ we have
$$
	d'_-f = -i \sum_{i=1}^m 
	S^{-1}(X_i) \triangleright f \otimes \sigma(K_\lambda)\gamma(X^i).
$$
\end{prop} 
\begin{pf}
The coproduct of
$S^{-1}(X_1)=-(1-q^{2nd_r}) K_\lambda^{-1} F_r$ 
is given by
$$
	S^{-1}(X_1) \otimes K_r^{-1} K_\lambda^{-1} +
	K_\lambda^{-1} \otimes S^{-1}(X_1). 
$$
Since $X_j=Y \triangleright X_1$ for some $Y \in U_q(\mathfrak{l})$ 
one obtains for $f \in \mathbb{C}_q[M]$ and 
$\sum_{i} g_i \otimes v_i \in \Gamma(M_q,\Sp)$ the relation
\begin{eqnarray}
&& \sum_{i} S^{-1}(X_j) \triangleright (g_i f) \otimes v_i
	\nonumber\\ 
&=& \sum_i (Y_{(3)} S^{-1}(X_1) S^{-1}(Y_{(2)}) \triangleright g_i)
	(Y_{(4)} K_r^{-1} K_\lambda^{-1} S^{-1}(Y_{(1)}) \triangleright f) 
	\otimes v_i \nonumber\\
&& + (Y_{(3)} K_\lambda^{-1} S^{-1}(Y_{(2)}) \triangleright g_i) 
	(Y_{(4)} S^{-1}(X_1) S^{-1}(Y_{(1)}) \triangleright f) 
	\otimes v_i\nonumber\\ 
&=& \sum_i (S^{-1}(X_j) \triangleright g_i) f 
	\otimes v_i + g_i (S^{-1}(X_j) \triangleright f) 
	\otimes K_\lambda \triangleright v_i,\nonumber 
\end{eqnarray} 
where we used the defining properties of 
$\mathbb{C}_q[M],\Gamma(M_q,\Sp)$ and the fact that 
$K_\lambda$ commutes with elements of $U_q(\mathfrak{l})$. 
\end{pf}
Since the right multiplication operators 
$R_g : f \mapsto fg$, $f,g \in \mathbb{C}_q[G]$ extend to
bounded operators on the Hilbert space obtained by completing 
$\mathbb{C}_q[G]$ with respect to Haar measure Proposition~\ref{df} implies: 
\begin{cor}\label{bounded}
The elements of $\Gamma_-'$ extend to bounded operators on $\H$.
\end{cor}
By [HK2, Corollary~5] there is a one-to-one correspondence between 
$m$-dimensional covariant differential calculi over $\mathbb{C}_q[M]$
and $m+1$-dimensional subspaces $\T$ of $\mathbb{C}_q[M]^\circ$ such
that
\begin{equation}\label{qts}
	\varepsilon \in \T,\quad
	\Delta(\T) \subset \T \otimes \mathbb{C}_q[M]^\circ,\quad
	U_q(\mathfrak{l}) \T \subset \T.
\end{equation}  
Here $\mathbb{C}_q[M]^\circ$ denotes the dual coalgebra of
$\mathbb{C}_q[M]$, see [HK2]. In view of [HK1, Theorem~6.5] it is 
sufficient to consider $\T \subset \pi(U_q(\mathfrak{g}))$, where
$\pi : \mathbb{C}_q[G]^\circ \rightarrow \mathbb{C}_q[M]^\circ$ is the
restriction map. Then 
$U_q(\mathfrak{l})\T \subset \T$ means that
$\pi(XY) \in \T$ for all $X \in U_q(\mathfrak{l})$
and $Y \in U_q(\mathfrak{g})$ with $\pi(Y) \in \T$.
The vector space $\T^0:=\{X \in \T \>|\> X(1)=0\}$ 
is called the quantum tangent space of the corresponding differential
calculus.
\begin{prop}
The vector space 
$\T_-^0 \subset \mathbb{C}_q[M]^\circ$ spanned
by $\pi \circ S^{-1}(X_i)$, $i=1,\ldots,m$ coincides with the quantum
tangent space of $\Gamma_-$.
\end{prop} 
\begin{pf}
For $f \in \mathbb{C}_q[M]$ we have
$$
	F_rK_rK_\lambda(f)=F_r((K_rK_\lambda) \triangleright f)=F_r(f)
$$
and similarly
$$
	(Y \triangleright X_1)(f)=
	YX_1(f)\quad\forall Y \in U_q(\mathfrak{l}).
$$
Hence the claim reduces to the fact that the tangent space of $\Gamma_-$ is
$U_q(\mathfrak{l}) \pi(F_r) \subset \mathbb{C}_q[M]^\circ$, see [HK1].
\end{pf}
It remains to show that $\Gamma_-'$ is
indeed isomorphic to $\Gamma_-$. To see this we first realize
$\Gamma_-$ as a calculus induced by a 
$m+1$-dimensional covariant differential calculus
over $\mathbb{C}_q[G]$. This calculus has tangent space
$$
	\T_-^{G0}:=\mathbb{C} S^{-1}(X_0) \oplus
	S^{-1}(\mathrm{ad}(U_q(\mathfrak{l}))X_1) 
	\subset U_q(\mathfrak{g}).
$$ 
Using that $S^{-1}$ is a coalgebra antihomomorphism
and that $K_\lambda$ commutes with all elements of
$U_q(\mathfrak{l})$ one calculates that
\begin{eqnarray}
&& \Delta(S^{-1}(Y \triangleright X_1)) \nonumber\\
&=& K_\lambda^{-1} \otimes S^{-1}(Y \triangleright X_1) +
	S^{-1}(Y_{(2)} \triangleright X_1) \otimes 
	S^{-1}(Y_{(1)}K_rK_\lambda S(Y_{(3)}))\nonumber\\
&\in& \T^G_- \otimes U_q(\mathfrak{g})\nonumber
\end{eqnarray}
for all $Y \in U_q(\mathfrak{l})$, where 
$\T^G_-:=\mathbb{C} \cdot 1 \oplus \T_-^{G0}$.
Therefore there is indeed a differential 
calculus $\Gamma^G_-$ over $\mathbb{C}_q[G]$ with quantum tangent
space $\T_-^{G0}$ (the last condition in (\ref{qts})
becomes trivial on quantum groups). By construction we have 
$\pi(\T_-^{G0})=\T_-^0$ and hence
$\Gamma^G_-$ induces $\Gamma_-$ [HK2, Corollary~9].\\
The general theory of 
covariant differential calculi over Hopf algebras with invertible
antipode (see [KS, Section~14.1]) implies that in $\Gamma^G_-$  
the differential can be written as
\begin{equation}\label{letzte}
	d^G_- f = \sum_{i=0}^m (S^{-1}(X_i) \triangleright f) \cdot
	\omega^i\quad
	\forall f \in \mathbb{C}_q[G],
\end{equation}    
where $\{\omega^i\}$ is a basis of $\Gamma^G_-$ consisting of invariant
1-forms. Proposition~\ref{df} 
generalizes the above formula to differential calculi over 
quantum flag manifolds.\\ 
The relation (\ref{letzte}) implies in particular that 
$$
	\sum_{i} f_i d^G_- g_i=0
	\quad \Leftrightarrow \quad 
	\sum_{i} f_i (S^{-1}(X_j) \triangleright g_i)=0\quad
	\forall j.
$$
The matrices 
$\sigma(K_\lambda)\gamma(X^i)$ are linearly
independent, because $\sigma(K_\lambda)$ is invertible,
$\gamma$ is injective and the $X^i$ are linearly independent. 
Furthermore $\mathbb{C}_q[G]$ is free
of zero divisors [J, 9.1.9]. Hence Proposition~\ref{df} implies
$$
	\sum_{i} f_i d_-' g_i = 0
	\quad \Leftrightarrow \quad 
	\sum_{i} f_i (S^{-1}(X_j) \triangleright g_i) = 0
$$ 
and we obtain:
\begin{prop}\label{fodc}
The map 
$$
	\psi : \Gamma_- \rightarrow \Gamma_-',\quad
	\sum_{j} f_j d_- g_j \mapsto \sum_{j} f_j d_-' g_j
$$
is an isomorphism of differential calculi.
\end{prop}

\end{document}